\documentclass[twoside,11pt]{amsart} \usepackage{epsfig}
\usepackage{latexsym} 
\usepackage{amsmath} \usepackage{amssymb}



\newtheorem{theorem}{Theorem}[section]
\newtheorem{lemma}[theorem]{Lemma}
\newtheorem{corollary}[theorem]{Corollary}
\newtheorem{remark}[theorem]{Remark}
\newtheorem{conjecture}[theorem]{Conjecture}

\newdimen\epsfxsize
\newdimen\epsfysize

\newcommand {\gap}     {\makebox[0.075 in]{}}   
\newcommand {\st}      {\gap : \gap}   

\newcommand {\set}[1]  {\left\{ {#1} \right\}}   
\newcommand {\ord}[1]  {{#1}^{\rm th}}

\newcommand{\C}{\mathbb{C}}

\newcommand{\Z}     {{\mathbb Z}}
\newcommand{\Zmod}  {\Z \bmod \Pi_k}
\newcommand{\Zmodp}  {\Z \bmod \Pi_{k+1}}
\newcommand{\N}[2]  {N_{{#2}}({#1})}
\newcommand{\E}[2]  {E_{{#2}}^{{#1}}}
\newcommand{\tC}[1] {{C^k_{#1}}}
\newcommand{\tE}[1] {{E^{k-1}_{#1}}}
\newcommand{\pgap}   {{\mathcal G}}
\newcommand{\gsum}[2]  {\Gamma_{{#1},{#2}}}

\begin{document}

\title{Expected gaps between prime numbers}

\date{2 Feb 06}

\author{Fred B. Holt}
\address{fbholt@u.washington.edu ;  4311-11th Ave NE \#500, Seattle, WA 98105}

\begin{abstract}
We study the gaps between consecutive prime numbers directly through
Eratosthenes sieve.  Using elementary methods, we identify
a recursive relation for these gaps and for  specific sequences of
consecutive gaps, known as constellations.  
Using this recursion we can estimate the numbers of a gap or of a constellation
that occur between a prime and its square.  This recursion also has explicit implications
for open questions about gaps between prime numbers, including three questions posed
by Erd\"os and Tur\'an.
\end{abstract}

\maketitle


\section{Introduction}

We work with the prime numbers in ascending order, denoting the
$\ord{k}$ prime by $p_k$.  Accompanying the sequence of primes
is the sequence of gaps between consecutive primes.
We denote the gap between $p_k$ and $p_{k+1}$ by
$g_k=p_{k+1}-p_k.$
These sequences begin
$$
\begin{array}{rrrrrrc}
p_1=2, & p_2=3, & p_3=5, & p_4=7, & p_5=11, & p_6=13, & \ldots\\
g_1=1, & g_2=2, & g_3=2, & g_4=4, & g_5=2, & g_6=4, & \ldots
\end{array}
$$

A number $d$ is the {\em difference} between prime numbers if there are
two prime numbers, $p$ and $q$, such that $q-p=d$.  There are already
many interesting results and open questions about differences between
prime numbers; a seminal and inspirational work about differences
between primes is Hardy and Littlewood's 1923 paper \cite{HL}.

A number $g$ is a {\em gap} between prime numbers if it is the difference
between consecutive primes; that is, $p=p_i$ and $q=p_{i+1}$ and
$q-p=g$.
Differences of length $2$ or $4$ are also gaps; so open questions
like the Twin Prime Conjecture, that there are an infinite number
of gaps $g_k=2$, can be formulated as questions about differences
as well.

A {\em constellation among primes} \cite{Riesel} is a sequence of consecutive gaps
between prime numbers.  Let $s=a_1 a_2 \cdots a_k$ be a sequence of $k$
numbers.  Then $s$ is a constellation among primes if there exists a sequence of
$k+1$ consecutive prime numbers $p_i p_{i+1} \cdots p_{i+k}$ such
that for each $j=1,\ldots,k$, we have the gap $p_{i+j}-p_{i+j-1}=a_j$.  
Equivalently,
$s$ is a constellation if for some $i$ and all $j=1,\ldots,k$,
$a_j=g_{i+j}$.

We will write the constellations without marking a separation between
single-digit gaps.  For example, a constellation of $24$ denotes
a gap of $g_k=2$ followed immediately by a gap $g_{k+1}=4$.
The number of gaps after $k$ iterations of the sieve is
$\Phi_k = \prod_{i=1}^k (p_i-1).$  For the small primes we will consider
explicitly, most of these gaps are single digits, and the separators introduce
a lot of visual clutter.  We use commas only to separate double-digit gaps in
the cycle.  For example, a constellation of $2,10,2$ denotes a gap of $2$
followed by a gap of $10$, followed by another gap of $2$.

We use elementary methods to study the gaps 
generated by Eratosthenes sieve directly.
By studying this sieve, we can 
estimate the occurrence of certain gaps and constellations
between $p_k$ and $p_k^2$.

From the methods developed below, we can calculate exactly
how many times a sequence $s$ of gaps occurs after $k$ stages of
Eratosthenes' sieve.  We don't know how many of these occurrences
will survive subsequent stages of the sieve to become constellations
among prime numbers.  However, for a prime $p$ we can make estimates for the
number that occur before $p^2$, all of which will survive as constellations among primes.
Thus our estimates and counts are only coincidentally commensurate with
tabulations against powers of ten.

The product of the first $k$ primes will be denoted by
$\Pi_k = \prod_{i=1}^k p_i.$

By the $p_k$-{\em sieve}, we mean those positive integers remaining
after removing all the multiples of the first $k$ prime numbers.
The $p_k$-sieve has a fundamental cycle of $\Phi_k$ elements modulo $\Pi_k$.
Most often we picture this fundamental cycle as the generators for $\Zmod$,
although it is also attractive to visualize 
these as the primitive $\Pi_k^{\rm th}$ roots of unity in $\C$.

\subsection{Organization of the material}
We proceed as follows.
We identify a recursive algorithm for producing each cycle of
gaps $\pgap(p_{k+1})$ from the preceding cycle $\pgap(p_k)$.
This recursion enables us to enumerate various gaps and constellations
in the $p_k$-sieve.
In the cycle of gaps $\pgap(p_k)$ of course,
all the gaps from $p_{k+1}$ and $p_{k+1}^2$ are actually gaps between
prime numbers.

We make a conjecture about the uniformity of the distribution of these
gaps and constellations.
From this conjecture we can make statistical estimates about the
expected number of occurrences of these gaps and constellations below $p_{k+1}^2$, 
and we compare these estimates with actual counts.

We make a weaker conjecture that every
constellation in $\pgap(p_k)$ occurs infinitely often as a
constellation among primes, provided the sum of the gaps in the constellation
is less than $2p_{k+1}$.
From this weaker conjecture we address several questions about gaps and differences
between prime numbers.  We show that Hardy and Littlewood's $k$-tuple conjecture 
on the differences between prime numbers \cite{HL} is equivalent to a conjecture on gaps.
We are also able to give exact answers to three questions posed by Erd\"os and Tur\'an \cite{ET}.

\subsection{New results}
This paper applies elementary methods to Eratosthenes sieve.  In a general sense, this is
well-trodden ground.  However, the specific insight of identifying the recursion on gaps appears
to be new.  We cast Eratosthenes sieve as a recursive operation directly on the cycle of gaps.
By studying this recursion we can enumerate particular gaps at every stage.  Moreover we
observe how much structure of the cycle of gaps at one stage of the sieve is preserved in subsequent
stages.  We can thereby easily enumerate the occurrences of specific constellations of primes that
have not previously been approachable (e.g. $2,10,2$).  

The conclusions on the recursion of gaps are precise.  To go further, we need to supplement
our rigorous work with an appropriate conjecture.  We first make a strong conjecture, that
under the recursion the copies of a specific constellation eventually approach a uniform distribution
in the cycle.  This conjecture on a uniform distribution allows us to make estimates of the occurrences
of constellations in the sieve as constellations among prime numbers.  These new estimates
compare favorably with existing estimates, and they allow us to estimate the occurrences of other
interesting constellations that have lain beyond the reach of existing techniques.

Backing off from the strong conjecture on uniformity, we make a weaker conjecture, that
sufficiently small constellations in the sieve occur infinitely often as constellations among
prime numbers.  This conjecture implies that the Twin Prime Conjecture is true.  However, it
goes further.  By identifying specific constellations and using the action of the recursion, we
answer three questions posed Erd\"os and Tur\'an \cite{ET}:
\begin{enumerate}
\item {\it Spikes.} $\limsup g_n/g_{n+1} = \infty$ and $\liminf g_n/g_{n+1}=0$
\item {\it Oscillation.} There is no $n_0$ such that for all $k \ge 1$,
 $g_{n_0+2k-1} < g_{n_0+2k}$ and $g_{n_0+2k} > g_{n_0+2k+1}$.
\item {\it Superlinearity.} $g_j < g_{j+1} < \ldots < g_{j+k}$ does have
infinitely many solutions for every $k$.
\end{enumerate}

The progress exhibited in this paper is the result of new elementary insights into
Eratosthenes sieve.  Specifically, we can track the cycle of gaps explicitly through
stages of the sieve.  While previous methods have jumped immediately to probabilistic
estimates, we examine the deterministic effect that the recursion has on subsequences
in the cycle of gaps.  Only after we have exhausted the exact results on these constellations, 
do we turn to simpler probabilistic estimates.

\section{Related Results}
There are of course several avenues of research into the distribution of
primes.  Research into constellations has been motivated primarily by
two conjectures:  the twin primes conjecture, and Hardy and Littlewood's
broader $k$-tuple conjecture \cite{HL,Rib}.

The twin primes conjecture asserts that the gap $g=2$ occurs infinitely often.
Work on this conjecture
has included computer-based enumerations \cite{IJ,NicelyTwins,PSZ}
and investigations of Brun's constant \cite{Sieves,Riesel,HL,Rib}.
Brun's constant is the sum of the reciprocals of twin primes.
This series is known to converge, 
and the sharpest current estimate \cite{Rib} is $1.902160577783278.$
One generalization of the twin primes conjecture is a conjecture by Polignac from
1849 \cite{Rib} that for every even positive integer $N$ 
there are an infinite number of gaps $g_k=N$.

Hardy and Littlewood \cite{HL,Rich} formulated their prime $k$-tuples
conjecture in this form:  if $b_1,\ldots,b_k$ is an admissible
$k$-tuple, then there are infinitely many $x$ such that
$x+b_1,\ldots,x+b_k$ are all prime. 
Admissibility in this context is
a condition on residue classes modulo smaller primes.

The work in \cite{HL} supports related conjectures estimating
the numbers of specific differences that should occur in the interval $[2,N]$
for any $N$.  For a difference $d$, their Conjecture B asserts that the number
$C_d(N)$ of prime pairs $(p,p+d)$ with $p \le N$ is asymptotically
\begin{eqnarray}\label{HLests}
C_d(N) 
 & \sim & 2 c_2 \prod_{q|d} \frac{q-1}{q-2} \int_2^N \frac{dx}{\ln^2 x}\\
 & \sim & 2 c_2 \frac{N}{\ln^2 N} \prod_{q|d} \frac{q-1}{q-2}
\end{eqnarray}
in which $q$ runs over the odd primes.  The constant $c_2$, known as the 
{\em twin prime constant} \cite{HL,Riesel} is given by the infinite product
$$c_2 = \prod \frac{p(p-2)}{(p-1)^2} = 0.6601618\ldots .$$
For an admissible $k$-tuple $b_1,\ldots,b_k$, Hardy and Littlewood \cite{HL}
conjectured the general estimate
\begin{eqnarray}\label{HLkests}
C_{b}(N) 
 & \sim & 2^k \prod_q \left(\frac{q}{q-1}\right)^{k+1}(1-\phi_q(b)/q)
 \; \cdot \; \int_2^N \frac{dx}{\ln^{k+1} x}
\end{eqnarray}
in which $q$ runs over the odd primes and
$\phi_q(b)$ is the number of distinct residue classes of
$0,b_1,\ldots,b_k$ modulo $q$.

Hardy and Littlewood's prime $k$-tuple conjecture addresses differences
among primes, but the primes in question need not be consecutive.
However the $k$-tuple conjecture has an equivalent
formulation as a conjecture on constellations, which we will see in 
Lemma~\ref{kLemma} below.
For the differences $d=2,4$ the estimates (\ref{HLests}) are also estimates
for the corresponding gaps $g=2,4$.

Most estimates for sequences of differences, e.g. \cite{HL,MV}, are derived by treating probabilities on residues as independent probabilities.  In contrast, the recursion identified below in Lemma \ref{RecursLemma} preserves the structure in the cycles of gaps at each stage of Eratosthenes sieve.  The occurrence of constellations in stages of the sieve is entirely deterministic.

Computational confirmation of these estimates has been carried out
by several researchers, notably in \cite{Brent3,NicelyTwins}.
The results in this paper and the tables and examples previously
published are not quite commensurate.
The tables and examples of \cite{Brent,IJ,Nicely,NicelyTwins} provide
estimates and counts of gaps with respect to large powers of ten.
Since we work directly with Eratosthenes sieve, our estimates
are given with respect to intervals $[p,p^2]$ for primes $p$.

Some researchers have applied their investigations of differences among
primes to study gaps and constellations.  The general surveys \cite{Rib,Riesel}
provide overviews of some of this work, and the estimate (\ref{HLkests})
from the seminal paper \cite{HL} can be used for constellations consisting of
$2$'s and $4$'s, e.g. $24$, $42$, $242$, $424$, etc.
Brent \cite{Brent} applied the principle of inclusion and exclusion
to the estimates (\ref{HLkests}) to obtain strong estimates for the
gaps $2,4,6,\ldots,80$.
Richards \cite{Rich} conjectured
that the constellation $24$ occurs infinitely often.
Clement \cite{quads} and Nicely \cite{NicelyQuads} have addressed the
constellation $242$, corresponding to {\em prime quadruplets}, pairs of twin
primes separated by a gap of $4$.

There are other lines of investigation \cite{Rib,GranICM} into the gaps between
prime numbers.
One line \cite{YP,Nicely,JLB,Brent2} has looked for 
the first occurrence of a gap.
Cramer \cite{Cramer} introduced probabilistic arguments to derive asymptotic
estimates for the sequence of gaps $g_k = \mathcal{O} (\ln^2 p_k).$
Quite recently, Green and Tao \cite{BGreen} have offered a proof that there
exist arbitrarily long sequences of primes in arithmetic progression.
By working with the convex hulls of the graphs of primes $(k,p_k)$ and of
the log-primes $(k,\ln p_k)$, Pomerance \cite{Pom} established a handful of
nice results about inequalities involving the arithmetic and geometric means of prime numbers.

\section{Recursion for the cycle of gaps}
The possible primes for the $3$-sieve are 
$$(1),5,7,11,13,17,19,23,25,29,31,35,37,41,43,\ldots$$
We investigate the structure of these sequences of possible primes
by studying the cycle of gaps in the fundamental cycle.
For example, $42$ is the cycle of gaps for the $3$-sieve.
We have
$$\pgap(3) = 42, \gap {\rm with} \gap g_{3,1}=4 \gap {\rm and} \gap 
g_{3,2}=2.$$

The lowest entry in the cycle of gaps is one less than the next prime:
$g_{k,1} = p_{k+1}-1$.
Denote the sum of the first $j$ gaps in $\pgap(p_k)$ by
$\gsum{k}{j}=\sum_{i=1}^j g_{k,i}$.
For the $p_k$-sieve, the $\ord{j}$ possible prime is given by 
$1+\gsum{k}{j}.$
Since we are studying the cycle of gaps in the $p_k$-sieve,
we know that there are $\Phi_k$ elements in one cycle, 
and the sum of the gaps in one cycle must be $\Pi_k$:
$$\gsum{k}{\Phi_k} = \Pi_k.$$

There is a nice recursion which produces $\pgap(p_{k+1})$ from
$\pgap(p_k)$.  We concatenate $p_{k+1}$ copies of $\pgap(p_k)$, and
add together certain gaps as indicated by the entry-wise product
$p_{k+1}*\pgap(p_k)$.  So the recursion consists of three steps.

\begin{lemma} \label{RecursLemma}
The cycle of gaps $\pgap(p_{k+1})$ is derived recursively from $\pgap(p_k)$.
Each stage in the recursion consists of the following three steps:
\begin{itemize}
\item[R1.] Determine the next prime, $p_{k+1} = g_{k,1} + 1$.
\item[R2.] Concatenate $p_{k+1}$ copies of $\pgap(p_k)$.
\item[R3.] Add together $g_{k,1}+g_{k,2}$, and record the index
for the location of this addition as $\tilde{i}_1=1$; 
for $n=1,\ldots,\Phi_k-1$, add $g_{k,j}+g_{k,j+1}$ and let 
$\tilde{i}_{n+1}=j$ if 
$$\gsum{k}{j}-\gsum{k}{\tilde{i}_n}= p_{k+1}*g_{k,n}.$$
\end{itemize}
\end{lemma}

\begin{proof}
We consider the cycle of gaps in relation to the generators of $\Z \bmod \Pi_k$.
Suppose the differences between 
consecutive generators in  $\Z \bmod \Pi_k$
is the cycle of gaps $\pgap(p_k)$.  By induction, we will show this relation holds
for $\pgap(p_{k+1})$ and $\Z \bmod \Pi_{k+1}$.

The next prime $p_{k+1}$ will be $1+g_{k,1}$, since this will be the smallest
integer both greater than $1$ and coprime to $\Pi_k$.

The second step of the recursion extends our list of possible primes
up to $\Pi_{k+1}+1$, the reach of the fundamental cycle for $p_{k+1}$.
For the gaps $g_{k,j}$ we extend the indexing on $j$ to cover these 
concatenated copies.  These $p_{k+1}$ concatenated copies of $\pgap(p_k)$
correspond to all the numbers from $1$ to $\Pi_{k+1}+1$ which are coprime
to $\Pi_k$.  For the set of generators of $\Pi_{k+1}$, we need only remove
the multiples of $p_{k+1}$. 

The third step removes the multiples of $p_{k+1}$.
Removing a possible prime amounts to
adding together the gaps on either side of this entry.
The only multiples of $p_{k+1}$ which remain in the copies of $\pgap(p_k)$
are those multiples all of whose prime factors are greater than $p_k$.
After $p_{k+1}$ itself, the next multiple to be removed will be $p_{k+1}^2$.

The multiples we seek to remove are given by $p_{k+1}$ times the generators
of $\Z \bmod \Pi_k$.  The consecutive differences between these will be given
by $p_{k+1} * g_{k,j}$, and the sequence $p_{k+1}*\pgap(p_k)$ suffices to cover
the concatenated copies of $\pgap(p_k)$.  We need not consider any fewer nor any
more multiples of $p_{k+1}$ to obtain the generators for $\pgap(p_{k+1})$.

In the statement of R3, the index $n$ moves through the copy of
$\pgap(p_k)$ being multiplied by $p_{k+1}$, and the indices $\tilde{i}_n$
mark the index $j$ at which the addition of gaps is to occur.
The multiples of $p_{k+1}$ in the $p_k$-sieve are given by
$p_{k+1}$ itself and $p_{k+1}*(1+\gsum{k}{j})$ for $j=1,\ldots,\Phi_k$.
The difference between successive multiples is $p_{k+1}*g_{k,j}$.
\end{proof}

%
\begin{figure}[t]
\begin{center}
\includegraphics[width=5in]{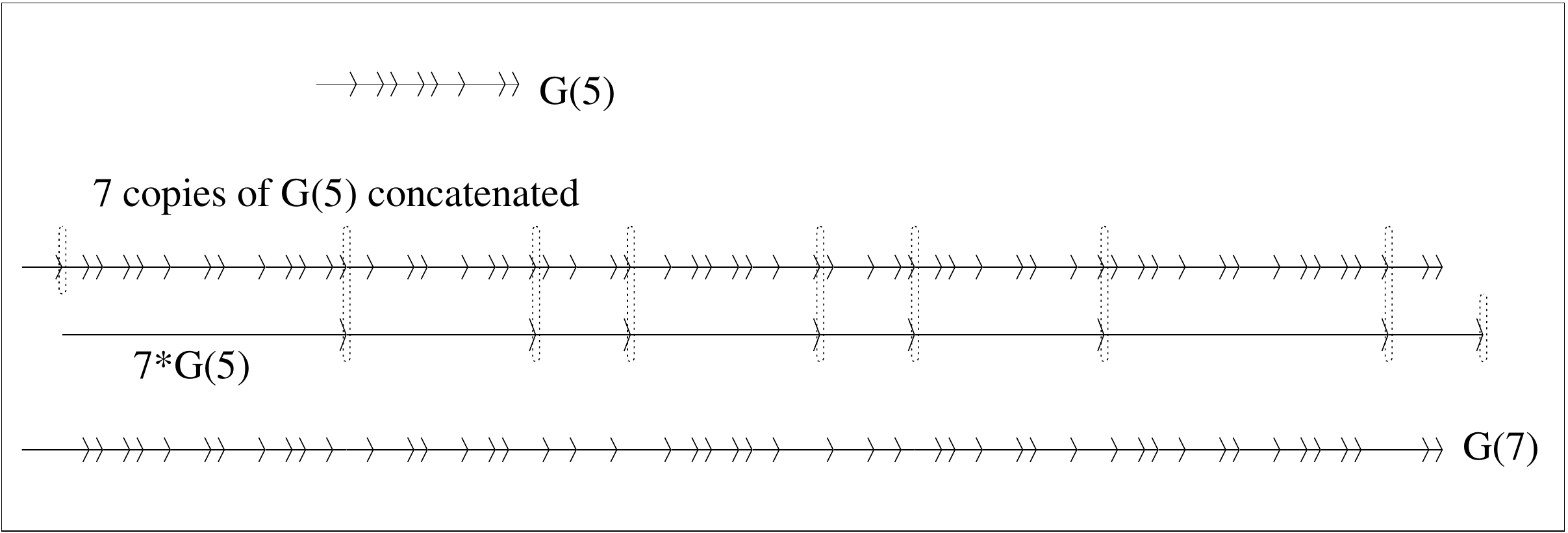}
\caption{\label{RecursFig} An illustration of the recursion that
produces the gaps for the next stage of Eratosthenes sieve.
The cycle of gaps $\pgap(7)$ is produced from $\pgap(5)$ by
concatenating $7$ copies, then adding the gaps indicated by
$7*\pgap(5)$.}
\label{default}
\end{center}
\end{figure}

\noindent {\bf Example: $\pgap(5)$.}
We start with $\pgap(3)=42.$
\begin{itemize}
\item[R1.]$p_{k+1}=5.$
\item[R2.] Concatenate five copies of $\pgap(3)$:
$$4242424242.$$
\item[R3.] Add together the gaps after $4$ and 
thereafter after cumulative differences of $5*\pgap(3)=20,10$:\\
\begin{tabular}{rcc}
$\pgap(5)$ & $=$ & $4+\overbrace{2424242}^{20}+\overbrace{42\;\;}^{10}$ \\
 & $=$ & $64242462.$ \\
\end{tabular} \\
Note that the last addition wraps over the end of the cycle
and recloses the gap after the first $4$.
\end{itemize}

\begin{remark}\label{EasyRmk}
The following results are easily established for $\pgap(p_k)$:

\begin{enumerate}
\item The first difference between additions is $p_{k+1}*(p_{k+1}-1)$,
which removes $p_{k+1}^2$ from the list of possible primes.
\item The last entry in $\pgap(p)$ is always $2$. 
This difference goes from $-1$ to $+1$ in $\Zmod$.
\item The last difference $p_{k+1}*2$ between additions, wraps from
$-p_{k+1}$ to $p_{k+1}$ in $\Zmodp$.
\item Except for the final $2$, the cycle of differences is symmetric:
$g_{k,j}=g_{k,\Phi_k-j}$.
\item If $g_{k,j}=g_{k,j+1}= \cdots = g_{k,j+m}=g$, then $g = 0 \bmod p$ for all
primes $p \leq m+2$.
\item The middle of the cycle $\pgap(p_k)$ is the sequence 
$$2^j,2^{j-1},\ldots,42424,\ldots,2^{j-1},2^j$$
in which $j$ is the smallest number such that $2^{j+1}>p_{k+1}$.
\end{enumerate}
\end{remark}

\noindent{\bf Example: $\pgap(7)$.} 
Following the steps in Lemma \ref{RecursLemma},
we construct $\pgap(7)$ from $\pgap(5)=64242462$.
This recursion is illustrated in Figure~\ref{RecursFig}.

\begin{itemize}
\item[R1.] Identify the next prime, $p_{k+1}= g_{k,1}+1 = 7.$
\item[R2.] Concatenate seven copies of $\pgap(5)$:
$$64242462 \; 64242462 \; 64242462 \; 64242462\; 64242462 \; 64242462 \;64242462$$
\item[R3.] Add together the gaps after the leading $6$ and 
thereafter after differences of $7*\pgap(5)$:
\begin{eqnarray*}
7*\pgap(5) &=& 42, 28, 14, 28, 14, 28, 42, 14 \\
\pgap(7) 
 &=&{\scriptstyle 
  6+\overbrace{\scriptstyle 424246264242}^{42}+
 \overbrace{\scriptstyle 4626424}^{28}+\overbrace{\scriptstyle 2462}^{14}+
 \overbrace{\scriptstyle 6424246}^{28}+\overbrace{\scriptstyle 2642}^{14}+
 \overbrace{\scriptstyle 4246264}^{28}+\overbrace{\scriptstyle 242462642424}^{42}+62} \\
 &=& {\scriptstyle 
 10, 2424626424 6 62642 6 46 8 42424 8 64 6 24626 6 4246264242, 10, 2}
\end{eqnarray*}
Note that the final difference of $14$ wraps around the end of the cycle,
 from the addition preceding the final $6$ to the 
addition after the first $6$.
\end{itemize}

\begin{theorem}\label{DelThm}
Each possible addition of adjacent gaps in the cycle $\pgap(p_k)$
occurs exactly once in the recursive construction of $\pgap(p_{k+1}).$
\end{theorem}

\begin{proof}
This is an implication of the Chinese Remainder Theorem.
Each entry in $\pgap(p_k)$ corresponds to one of the generators
of $\Z \bmod \Pi_k$. The first gap $g_{k,1}$ corresponds to $p_{k+1}$, and
thereafter $g_{k,j}$ corresponds to $1+\gsum{k}{j}$.
These correspond in turn to unique combinations of nonzero
residues modulo the primes $2,3,\ldots,p_k$.  In the $p_{k+1}$
copies of $\pgap(p_k)$, each copy of a particular gap $g_{k,j}$
has its combination of residues augmented by a unique residue modulo
$p_{k+1}$.  Exactly one of these has residue $0 \bmod p_{k+1}$, so we
perform $g_{k,j}+g_{k,j+1}$ for this copy and only this copy of $g_{k,j}$.
\end{proof}

\begin{corollary}
In $\pgap(p_{k+1})$ there are at least two entries of $2p_k$.
\end{corollary}

\begin{proof}
In forming $\pgap(p_k)$, we concatenate $p_k$ copies of
$\pgap(p_{k-1})$.  At the transition between copies we have the
subsequence $(p_k-1) 2 (p_k-1)$.  In $\pgap(p_k)$ each of
the two additions takes place, so the sequences $(p_k-1)(p_k+1)$ and
$(p_k+1)(p_k-1)$ both occur.  In $\pgap(p_{k+1})$ the addition
in each of these two sequences occurs in one of the $p_{k+1}$
copies.
\end{proof}

This corollary provides long runs of composite numbers earlier than the
traditional elementary constructions. For example, suppose we were looking
for runs of one thousand consecutive composite numbers.  
In the traditional approach, we would take $p_{169} = 1009$
and note that
 $$\{\Pi_{169} + 2, \Pi_{169} +3, \ldots, \Pi_{169}+1009, \Pi_{169}+1010 \}$$
and
 $$\{\Pi_{169} -1010, \Pi_{169} -1009, \ldots, \Pi_{169}-3, \Pi_{169}-2 \}$$
are runs of $1009$ composite numbers.
Using the above corollary, we take $p_{96} = 503$ and note that in
$\{p_{97}, \ldots, \Pi_{97}\}$ there occur at least two runs 
of $1006$ composite numbers.

\section{Specific constellations in $\pgap(p)$}.

A {\em constellation} is a sequence of gaps.
In this section, we use the recursion of Lemma \ref{RecursLemma}
to count the number of occurrences of a constellation
in $\pgap(p_k)$.  Then in the next section we estimate how
many of these occurrences survive as constellations among
primes smaller than $p_{k+1}^2$.

Whether a constellation $s$ continues to occur in $\pgap(p)$ under
the recursion depends on the number of gaps.

\begin{lemma}\label{SurviveLemma}
Let $s$ be a constellation of $j$ gaps in $\pgap(p_k)$.
If $j < p_{k+1}-1$, then copies of $s$ will appear in all
$\pgap(p)$ with $p \ge p_k$.
\end{lemma}

\begin{proof}
By Lemma~\ref{RecursLemma} any constellation is initially replicated
$p_{k+1}$ times, in Step R2.  But then by Theorem~\ref{DelThm}, each
possible addition occurs exactly once, corrupting up to $j+1$ copies
of $s$.  If $j < p_{k+1}-1$, then at least one copy of $s$ survives
intact.
\end{proof}

By Lemma~\ref{SurviveLemma}, if a constellation is short enough, then it 
will propagate through the $\pgap(p)$ until the conditions of the
following theorem are met, and from that point on we can enumerate the
occurrences of the constellation by the recursive equations provided
in the theorem.

\begin{theorem}\label{CountThm}
Let $s$ be a constellation of $j$ gaps in $\pgap(p_k)$, such that 
the sum of these $j$ gaps is less than $2p_{k+1}$.
Let $S$ be the set of all constellations $\bar{s}$ which would produce
$s$ upon one addition of differences.
Then the number $\N{p}{s}$ of occurrences of $s$ in $\pgap(p)$
 satisfies the recurrence
$$\N{p_{k+1}}{s} = (p_{k+1}-(j+1))*\N{p_k}{s}
 + \sum_{\bar{s} \in S} \N{p_k}{\bar{s}}$$
\end{theorem}

\begin{proof}
We account for the number of copies of $s$ which
survive the recursion intact, and we add to this 
the number of new copies of $s$ generated from other
sequences.

The $j$ gaps in $s$ can be closed in $j+1$ ways.
Theorem \ref{DelThm} tells us that each of these will occur
exactly once.
If the sum of the gaps in $s$ is less than $2p_{k+1}$,
then these closings are guaranteed to occur in distinct copies of $s$.
This establishes the first term on the right-hand side.

Finally, the summation on the right-hand side
accounts for occurrences which are generated from
other constellations.
\end{proof}

\begin{figure}[tb] 
\centering
\includegraphics[width=5.25in]{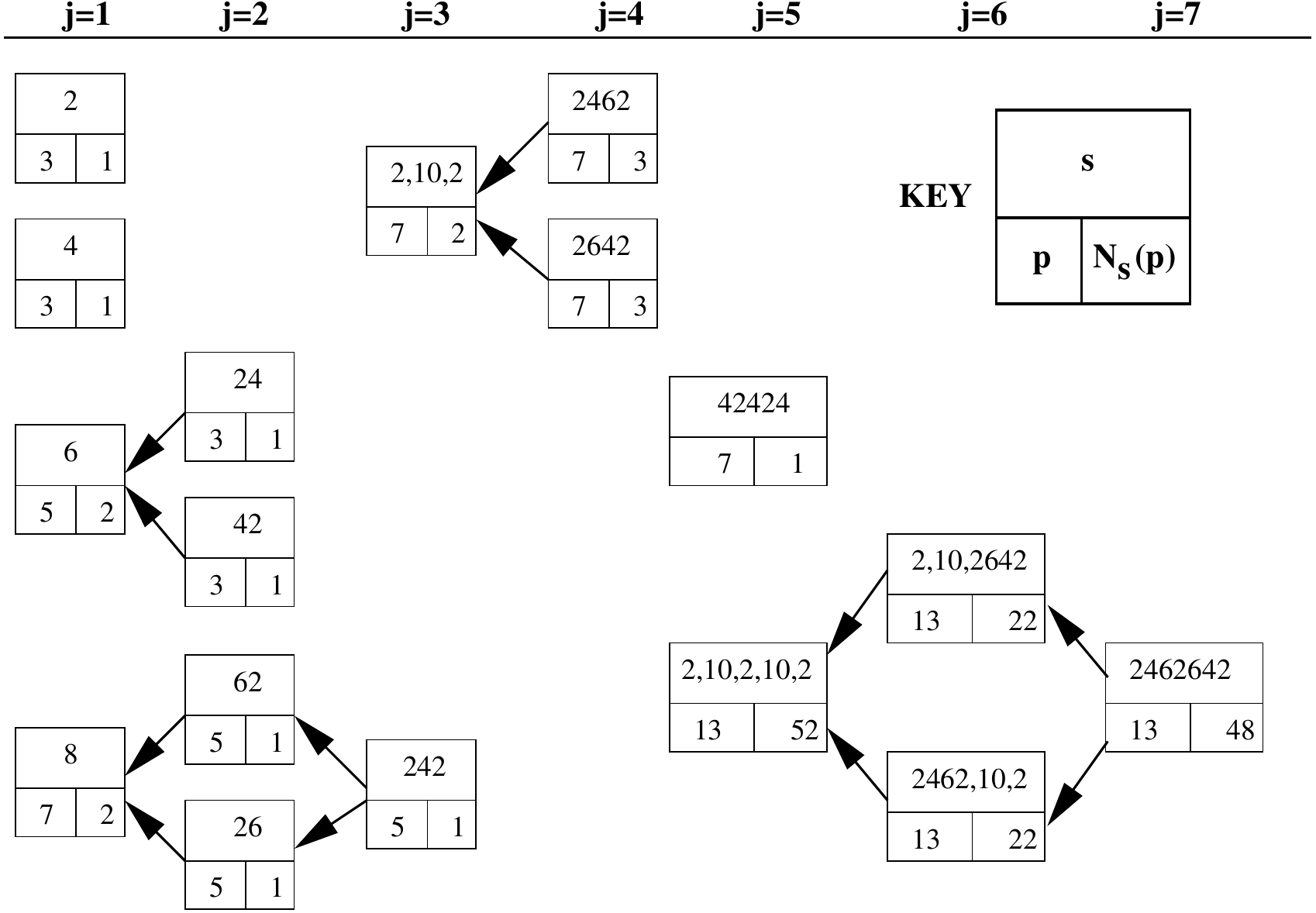}
\caption{\label{GrowthFig} The numbers $\N{p}{s}$ of a constellation $s$ 
are completely described by figures like this.  All constellations of $j$
gaps have identical dominant terms in Theorem~\ref{CountThm}, so we line
up our entries for constellations in columns indexed by $j$.  For each
constellation $s$ we include its initial conditions:  the first prime $p$
for which $s$ occurs in $\pgap(p)$ and the conditions of the theorem hold;
 and the number of occurrences of $s$ in $\pgap(p)$.
We indicate driving terms for the recurrence in Theorem~\ref{CountThm}
by directed arrows. Despite any disparities in initial conditions between
two constellations of $j$ gaps, the constellation with more driving terms
will rapidly become more numerous.}
\end{figure}

\begin{corollary}\label{NjgapCor}
Let $s$ be a constellation of $j$ gaps in $\pgap(p_k)$.  
If $j < p_{k+1}-1$, then the number of copies of $s$, $\N{p}{s}$,
will be dominated by the factors $\prod(p-j-1)$.
Thus all constellations of $j$ gaps grow asymptotically at
the same rate.
\end{corollary}

Although the asymptotic growth rates of all constellations of $j$ gaps
are equal, the initial conditions and driving terms are important.
Brent \cite{Brent} made analogous observations for single gaps ($j=1$).
His Table 2 indicates the importance of the lower-order effects in
estimating relative occurrences of certain gaps.

Figure~\ref{GrowthFig} provides a diagram of growths among various
constellations.  For each constellation $s$, its growth is dominated
by the number $j$ of gaps through the factors $(p-j-1)$.  
Those constellations $S$ which produce $s$ after one addition
provide the driving terms in the summation, but they are added in without
a multiplier, and they themselves grow with factors $(p-j-2)$, since
they are one gap longer than $s$.  Finally, the counts $\N{p}{s}$ for
various constellations will vary by their initial conditions.  These
initial conditions are the prime $p$ for which $s$ first occurs in $\pgap(p)$
and the number of these first occurrences.

Neither $2$'s nor $4$'s are generated from any other sequence.
In $\pgap(3)$
there is one $2$ and one $4$, so the numbers of occurrences
of these two gaps will continue to be equal through all stages of the
sieve.

From Figure~\ref{GrowthFig} or from the symmetry of $\pgap(p)$,
we know that the constellations $24$ and $42$ will always occur 
equally often.  However, $242$ has no driving terms and so it will
soon be outnumbered by the constellation $2,10,2$.
The constellation $242$ also generates occurrences of the constellations
$26$ and $62$, and these in turn generate $8$'s; a gap of $6$ is generated
by the constellations $24$ and $42$, which themselves have no generators.
Under the recursion we will eventually have
$$\N{p}{2} = \N{p}{4} < \N{p}{6} < \N{p}{8}.$$

More surprising constellations include $42424$ or $2,10,2$ or even
$2,10,2,10,2$.
Due to the terms for the sequences $2462$ and $2642$, 
the constellation $2,10,2$ becomes more abundant than the constellation
$242$.  These two equal additional terms in the recursion for $2,10,2$
fade in their significance compared to the multiplier $(p-4)$.  By $p=89$,
these two terms are contributing around one-half a percent of the total.
Similarly, while the constellations $42424$ and $2,10,2,10,2$ both occur,
ultimately $2,10,2,10,2$ will be more abundant due to driving terms.

\begin{center}
\begin{tabular}{r|rrr|rrr|r} \hline
$p$ & $\N{p}{2}$ & $\N{p}{24}$ & $\N{p}{6}$ 
 & $\N{p}{242}$ & $\N{p}{26}$ & $\N{p}{8}$ & $\N{p}{2,10,2}$ \\ \hline
$ 5$ &    $ 3$ &    $2$ &     $2$ &    $1$ &    $1$ &    $0$ &    $0$ \\
$ 7$ &    $15$ &    $8$ &    $14$ &    $3$ &    $5$ &    $2$ &    $2$ \\
$11$ &   $135$ &   $64$ &   $142$ &   $21$ &   $43$ &   $28$ &   $20$ \\
$13$ &  $1485$ &  $640$ &  $1690$ &  $189$ &  $451$ &  $394$ &  $216$ \\
$17$ & $22275$ & $8960$ & $26630$ & $2457$ & $6503$ & $6812$ & $3096$ \\ \hline
\end{tabular}
\end{center}

In the calculation for $\N{p}{2,10,2}$, we can also observe the necessity
of the requirement that the sum of the elements in a constellation be
less than $2p$.  One initial condition for $\N{p}{2,10,2}$ is $\N{7}{2462}=3$.
Since there are four elements in the constellation $2462$, we expect the
number of these to grow as $(p-5)$, and without checking, we might assume
that there would be $2$ such constellations in $\pgap(7)$.  However, the
sum of these elements is $14=2*7$, and in the construction of $\pgap(7)$
a sum of $2*7$ falls perfectly across one copy of $2462$.  The result is
that two of the gaps are closed in this single copy of $2462$, letting an
additional one of the seven original copies (from step R2) survive.  Until
$p$ is greater than twice the sum of the elements in the constellation,
we cannot be certain that each of the gaps in this constellation will be
closed in distinct copies.

From the work above, we have methods for tracking the exact number of copies
of any constellation through the stages of Eratosthenes sieve.  Unless the
constellation occurs in $\pgap(5)$, it must initially result from the closing of gaps
in longer constellations.  When a constellation appears in some $\pgap(p)$,
Lemma~\ref{SurviveLemma} and Theorem~\ref{CountThm} provide the conditions
and formulae for the exact number of this constellation that will appear in each 
stage of the sieve.  This system is illustrated in Figure \ref{GrowthFig}.

\section{Expected Constellations of Prime Numbers}

All of the work above is deterministic.  We have identified a recursion that
produces the cycle $\pgap(p_k)$ of gaps between consecutive generators of 
$\Z \bmod \Pi_k$.  We are able to identify and
count individual gaps and constellations as the recursion progresses.
In this section we introduce two conjectures, to help us explore the potential of
this approach.

The quantities $\N{p}{s}$ provide exact counts of the constellation $s$ in the cycles $\pgap(p)$
as we iterate on $p$.  These counts indicate how often $s$ occurs as a constellation among
possible primes.  We would like to use this information to draw conclusions about how often
$s$ occurs as a constellation among prime numbers.

To estimate how often $s$ occurs as a constellation among prime numbers, we combine two insights.
First, all constellations in $\pgap(p_k)$ from $p_{k+1}$ to $p_{k+1}^2$ are constellations among
primes.  Second, the recursion suggests uniform distributions.

\begin{remark}
In $\pgap(p_k)$, all the gaps that occur after $g_{k,1}=p_{k+1}-1$ and 
before $p_{k+1}^2$ are
actually gaps between prime numbers. 
For $j>1$,
if $ \gsum{k}{j} < p_{k+1}^2$, then $g_{k,j} = g_{k+j}.$
\end{remark}

Constellations which fall in the interval $[p,p^2]$ are constellations among primes.
We want to understand
the distribution of the $\N{p}{s}$ copies of a constellation $s$ over the cycle from $1$ to $\Pi_k$.
Observe that the concatenation step in the recursion produces
uniformly distributed copies of every constellation.  By using the interval $[p, p^2]$ to
sample these approximate uniform distributions of constellations, we can make estimates $E_s(p)$
of how often $s$ occurs as a constellation among primes.

\subsection{Nearly Uniform Distributions}
The recursion in Lemma~\ref{RecursLemma} suggests that if we track the images
of some gap $g$ through several stages of Eratosthenes sieve, these images
will be almost uniformly distributed in the fundamental cycle.
In $\pgap(p_k)$ pick any gap $g=g_{k,j}$.
Step R2 of the recursion creates $p_{k+1}$ copies of this gap, uniformly
distributed in the interval $[1,\Pi_{k+1}]$.  Step R3 removes two of these
copies of $g$.  In $\pgap(p_{k+2})$, step R2 creates $p_{k+2}$ copies
of the set of nearly uniformly distributed $p_{k+1}-2$ copies of $g$.  These
$p_{k+2}$ copies are uniformly distributed in the interval $[1,\Pi_{k+2}]$.
After step R3, the $(p_{k+2}-2)(p_{k+1}-2)$ images of $g$ in $\pgap(p_{k+2})$
are approximately uniformly distributed.  As we continue applying the
recursion, the abundant images of $g$ are pushed toward uniformity by
step R2 and trimmed symmetrically by step R3.
In Figure \ref{RecursFig} we can see this effect after only one stage of the
sieve.

The preceding observations apply to constellations as well.  For a 
constellation, the distribution of its images may not well-approximate a
uniform distribution until several of the recursions have been applied.
If the constellation contains $j$ gaps, then after step R2 in the recursion
distributes $p$ copies of $s$ uniformly, step R3 trims $j+1$ of these
copies by closing gaps in the constellation.  These $j+1$ copies
of $s$ are trimmed in symmetrical fashion.  Although $j$ may initially
be almost as large as $p$, through the recursion $p$ grows while $j$
remains fixed.  Before long, the trimming of step R3 will not substantially
disrupt the uniformity of the distribution enforced by step R2.

These observations, about the effects of the recursion on the distribution
of occurrences of a particular constellation, indicate a uniform distribution, but
there is more work to be done in understanding the distribution of copies of $s$
in $\pgap(p)$.  Further investigation will proceed along two lines:  computer searches
and statistical analyses.  For our present purposes we leave the
desired result as a conjecture.

\begin{conjecture}\label{uniconj}
Under the recursion in Lemma \ref{RecursLemma}, all constellations
in $\pgap(p)$ of sum less than $2p$ tend toward
a uniform distribution in $\pgap(P)$ for all primes 
$P \gg p$.
\end{conjecture}

Our estimates of the frequency of occurrences of certain constellations
require the uniform distribution asserted by this conjecture.  However,
other results, for example on some conjectures by Erd\"os and Tur\'an \cite{ET},
require a much weaker conjecture.

\begin{conjecture}\label{surviveconj}
Under the recursion in Lemma \ref{RecursLemma}, all constellations
in $\pgap(p_k)$ of sum less than $2p_{k+1}$ and with fewer than $p_{k+1}-1$ gaps 
occur infinitely often as
constellations among larger primes.
\end{conjecture}

From Lemma \ref{SurviveLemma}, we know that copies of these constellations will
survive throughout the $\pgap(p)$, and from Theorem \ref{CountThm} that the
number of copies will grow superexponentially.  Our first conjecture postulates
a distribution which enables us to count the copies which fall in
the intervals $[p,p^2]$ for every $p$.  

This first conjecture suggests that the third step in the recursion has an approximately uniform effect on the uniform distribution of copies created by the second step.
The second weaker conjecture asserts only that for some subsequence
of the primes, each interval $[p,p^2]$ contains at least one
copy of the constellation.  Suppose the third step in the recursion removes more copies of a constellation from the ends of the cycle of gaps, so that the superexponential number of copies are clustered in the middle (recall that $\pgap(p)$ is symmetric).  This second conjecture postulates that occasionally a copy of the constellation falls into the interval $[p,p^2]$. 

To bolster these conjectures, we can step back and look at the aggregate population of all constellations.  Not all constellations can accumulate in the middle of $\pgap(p)$.  Some must fall
in $[p,p^2]$.  If the distribution for some particular constellation is forever biased strongly toward
the middle of the cycle, then the distributions of some other constellations must compensate for
this bias.  That is, if some constellations fall below the expected number of occurrences, then 
other constellations must exceed the expectations.

\vskip .25in

\subsection{Estimates}\label{EstSection}
For a given constellation $s$, we use the approximate uniformity
of the $\N{p_k}{s}$ copies of $s$ in $\pgap(p_k)$ to estimate the number
that survive as constellations among primes.  There are $\Phi_k$ gaps
in $\pgap(p_k)$.  The average length $\mu$ of a gap in $\pgap(p_k)$ is
$$ \mu = \frac{\Pi_k}{\Phi_k} \sim e^\gamma \ln p.$$
The limit is due to Mertens, and its derivation is recorded in \cite{HW}.
In an interval $I$ in $[1,\Pi_k]$, we can expect there to be 
$\left|I\right|/\mu$ gaps.
Thus the expected number of constellations $s$ in an interval $I$ is
given by
$$ E(s,I) = \frac{\N{p_k}{s}}{\Phi_k - j+1} \frac{\left|I\right|}{\mu}.$$
The first factor is the fraction of gaps in $\pgap(p_k)$ that start
a copy of $s$, and the second factor is the expected number of gaps
in the interval $I$.  As $k$ gets large, the constant correction of $j-1$ 
in the first denominator becomes inconsequential.

Of particular interest to us is the interval $[p_{k+1},p_{k+1}^2]$:
\begin{equation}\label{EsEq}
\E{k}{s} = E(s,[p_{k+1},p_{k+1}^2]) 
\sim \frac{\N{p_k}{s}}{\Phi_k} 
 (p_{k+1}^2-p_{k+1}) \frac{e^{-\gamma}}{\ln p_k}.
\end{equation}

As a first application of this approach, we estimate the number of
twin primes between $p$ and $p^2$.
\begin{eqnarray*}
\E{k}{2} &=& (p_{k+1}^2-p_{k+1}) 
\frac{1}{\mu} \frac{\prod (p_j -2)}{\prod (p_j-1)} \\
 & \sim & 2 e^{-2\gamma}c_2 \frac{p_{k+1}^2-p_{k+1}}{\ln^2 p_k}
\end{eqnarray*}

Our estimate for the number of twin primes must be contrasted with
those estimates provided by Hardy and Littlewood, which are supported by 
vast computation \cite{Brent,NicelyTwins}:
$$
\begin{array}{lcc}
\# \set{g_i = 2 \st p_i \in [p_{k+1},p_{k+1}^2]} & 
 \sim & 2 e^{-2\gamma}c_2 \frac{p_{k+1}^2-p_{k+1}}{(\ln p_k)^2} \\
{\rm vs.} & & \\
\# \set{g_i = 2 \st p_i \in [2,N]} & 
 \sim & 2 c_2 \int_2^N \frac{dx}{\ln x} 
 \; \sim \;  2 c_2 \frac{N}{(\ln N)^2}
\end{array}
$$

Some actual counts $C^k_2$ of the single-gap constellation $2$ occurring between
$p_k$ and $p_k^2$ are tabled below.
These actual counts are compared to the estimates from the preceding
stage of the sieve, $\E{k-1}{2}$, and to the Hardy-Littlewood estimate $HL^k_2$.
The right half of the table compares the counts and estimates for the single-gap
constellations $6$ and $8$.

\begin{center}
\begin{tabular}{r|rrr||rr|rr} \hline
$p_k$ & $\tC{2}$ & $\tE{2}$  & $HL^k_2$ & $\tC{6}$ & $\tE{6}$ &  $\tC{8}$ & $\tE{8}$ \\ \hline
$11$ &  $8$ &     $8$ &  $4$ & $7$ & $7$ & $2$ & $1$ \\
$13$ &  $9$ &     $9$ &  $6$ & $10$ & $10$ & $1$ & $2$ \\
$101$ & $202$ &  $181$ & $152$ & $296$ & $286$ & $104$ & $96$ \\
$199$ & $574$ & $530$ & $457$ & $898$ & $878$ & $335$ & $312$ \\
$499$ & $2557$ & $2470$ & $2112$ & $4099$ & $4263$ & $1672$ & $1579$ \\
$1009$ & $8278$ & $8217$ & $6997$ &  $13715$ & $14521$ & $5643$ & $5506$ \\
 $1999$ & $26777$ & $26742$ & $22788$ & $44785$ & $48159$ & $18762$ & $18601$ \\ 
$2503$ & $39326$ & $39558$ & $33717$ & $66333$ & $71628$ & $27924$ & $27811$ \\
$4999$ & $130343$ & $133426$ & $113623$ & $223691$ & $245166$ & $96283$ & $96528$ \\
$10007$ & $440666$& $457406$ & $389427$ & $769389$ & $850965$ & $334491$ & $338959$ \\
 $12503$ & $653634$ & $681311$ & $579620$ & $1146148$ & $1271986$ & $499702$ & $508315$\\
 $14939$ & $895790$ & $936917$ & $797157$ & $1576337$ & $1753990$ & $689398$ & $702709$ \\
\hline
\end{tabular}
\end{center}

Twin primes occur in interesting constellations, for example
the prime quadruplets constellation $s=242$.
This constellation occurs in $\pgap(5)$.
With $\N{5}{242}=1$, we calculate
the number of expected occurrences of the
constellation $242$ between $p$ and $p^2$.
Under the recursion, $\N{p_{k+1}}{242} = (p_{k+1}-4)\N{p_k}{242}$.
So, under our conjecture of uniformity, the expected number of these
constellations between $p_{k+1}$ and $p_{k+1}^2$ is
\begin{eqnarray} \label{E242}
\E{k}{242} &=& \frac{(p_{k+1}^2-p_{k+1})}{\mu} 
 \prod_{p_j=5}^{p_k} (p_j-4) / \Phi_k \\ \nonumber
 & = & \frac{(p_{k+1}^2-p_{k+1})}{\mu^4} \prod_{p_j=5}^{p_k} (p_j-4) \cdot 
\frac{\Pi_k^3}{\Phi_k^4} \\ \nonumber
 & \sim & \frac{27}{2} c_4 e^{-4\gamma}\frac{(p_{k+1}^2-p_{k+1})}{\ln^4 p_k},\\ 
\nonumber
{\rm with} \; c_4 &=& \prod_{q\ge 5} \frac{q^3(q-4)}{(q-1)^4} = 0.30749\ldots 
\end{eqnarray}
This constant $c_4$ can be derived from the estimates (\ref{HLests}) 
and (\ref{HLkests}). See for example \cite{Riesel}.

The extensive computations of \cite{NicelyQuads} support the Hardy-Littlewood
estimates \cite{Riesel} for $s=242$:
$$ \frac{27}{2} c_4 \int_2^N \frac{dx}{\ln^4 x} \; \sim \; \frac{27}{2}c_4 \frac{N}{(\ln N)^4}.$$

Calculating $\E{k}{2,10,2}$ is more involved because of the driving
terms from the constellations $2462$ and $2642$.  We note that
for all $p\ge 13$, $\N{p}{2,10,2} > \N{p}{242}$, but that both have the
same dominating factor of $p-4$. 

\begin{center}
\begin{tabular}{r|rrr||rr|rr} \hline
$p_k$ & $\tC{242}$ & $\tE{242}$  & $HL^k_{242}$ & $\tC{2,10,2}$ & $\tE{2,10,2}$  & 
$\tC{2,10,2,10,2}$ & $\tE{2,10,2,10,2}$  \\ \hline
$11$ &  $2$ & $2$ &  $0$ & $0$ & $0$ & $0$ & $0$ \\
$13$ &  $1$ & $1$ &  $0$ & $1$ & $1$ & $0$ & $0$ \\
$101$ & $10$ & $9$ & $5$ & $18$ & $16$ & $1$ & $1$ \\
$199$ & $20$ & $20$ & $12$ & $35$ & $37$ & $2$ & $2$ \\
$499$ & $56$ & $67$ & $42$ & $118$ & $135$ & $5$ & $6$ \\
$1009$ & $167$ & $182$ & $114$ & $325$ & $377$ & $10$ & $13$ \\
 $1999$ & $459$ & $490$ & $308$ & $873$ & $1041$ & $25$ & $29$ \\ 
 $2503$ & $620$ & $683$ & $431$ & $1249$ & $1464$ & $38$ & $39$ \\
 $4999$ & $1714$ & $1948$ & $1228$ & $3621$ & $4255$ & $84$ & $95$\\
 $10007$ & $4760$& $5712$ & $3604$ & $10502$ & $12686$ & $212$ & $243$\\
 $12503$ & $6657$ & $8118$ & $5114$ & $14872$ & $18113$ & $300$ & $331$ \\
 $14939$ & $8777$ & $10753$ & $6777$ & $19556$ & $24079$ & $378$ & $424$\\
 \hline
\end{tabular}
\end{center}

The tables above sample the output from a program that searches for prime constellations,
up to the limits of long integer arithmetic and available RAM.   The largest prime for which 
constellations can be counted accurately by this program is $14939$.
Extending these computational results is one open avenue of research.

Correlations among copies of a constellation are somewhat preserved by the recursion.
Our conjecture~\ref{uniconj} does not take into account these correlations.
For example, the constellation $2,10,2,10,2$ occurs in $\pgap(7)$ and thereafter.
This constellation contains two occurrences of $2,10,2$ and
will consequently introduce jumps in $C_k(2,10,2)$. 
Any deviations from a uniform distribution are altered during the recursions.
These changes to the distribution are therefore occurring on the scale of $\Pi_k$.
As a reference point for this scale, the results tabulated in this paper, through $p=14939$
require only $\Pi_9$, with $p_9 = 23$.

\subsection{Refinements}
Based on our conjecture about uniformity, we could calculate the
expected number of a constellation $s$ that fall between
$p_{k+1}$ and $p_{k+1}^2$ as
$$ \E{k}{s} = (p_{k+1}^2-p_{k+1})\frac{\N{p_k}{s}}{\Pi_k}.$$

For small primes we can improve these estimates.
The symmetry of $\pgap(p)$ allows us to refine the denominator 
to $\Pi_k - 2p_{k+1}$.  Observations about the middle of $\pgap(p)$
allow a further refinement in the denominator.  
Our estimates do not need to use the entire interval $[1, \Pi_k]$.
The cycle $\pgap(p)$ is symmetric, and its middle constellation
is $2^j \ldots 4 2 4 2 4 \ldots 2^j$.  Thus we could
work over the interval $[p_k,\Pi_k/2-2^{j+1}]$ instead;
the term $2^{j+1}$ is the smallest power of $2$ greater than $p_{k+1}$.

These refinements may
require us to decrement the numerator slightly, adjusting for occurrences
we know we've excluded, e.g. the final $2$.
We could adjust the interval $[p,p^2]$ to $[p,p^2-l]$ to adjust for the sum $l$
of the gaps in the constellation.
These refinements are constant or of order $p$, which become
inconsequential rapidly in the face of factors like $\Pi_k$.

\section{Implications for other problems on constellations}\label{ImpSection}
Although Hardy and Littlewood's prime $k$-tuple conjecture 
is formulated about differences between
primes, we show here that the conjecture 
has an equivalent
formulation as a conjecture on constellations.

To the $k$-tuple we associate a constellation.
Without loss of generality we can assume the $k$-tuple is in
ascending order.  From this $k$-tuple we can derive
a sequence of $k-1$ differences: 
\begin{equation}\label{kconst}
s_b = (b_2-b_1),(b_3-b_2),\ldots,(b_k-b_{k-1}).
\end{equation}
This sequence is the constellation we want to associate with the $k$-tuple.
Notice, however, that the $k$-tuple conjecture does not require
that the primes be consecutive.  So while the $k$-tuple conjecture
can be trivially reformulated for a sequence of differences, we
have a little work to do to establish that this is equivalent to a
conjecture on constellations.

\begin{lemma}\label{kLemma}
The $k$-tuple conjecture is equivalent to the conjecture
that the constellation $s_b$ occurs in $\pgap(p)$ for
some $p$ with $2p > b_k-b_1$.
\end{lemma}

\begin{proof}
For the $k$-tuple $B=(b_1,\ldots,b_k)$ we put the $b_i$ in ascending order,
and then associate this with the constellation
$$s_b = (b_2-b_1),(b_3-b_2),\ldots,(b_k-b_{k-1}).$$

$(\Longleftarrow)$ If the constellation $s_b$ occurs in $\pgap(p)$ with
$2p > b_k-b_1$, then this constellation occurs infinitely often as
a constellation among primes, and the $k$-tuple conjecture would be true
for $B$.

$(\Longrightarrow)$ Conversely, suppose the $k$-tuple conjecture were true
for $B$.
Then the sequence of
differences in (\ref{kconst}) occurs infinitely often among the primes.
Let $p$ be a prime large enough so that two conditions hold:
$2p > b_k-b_1,$ and
one instance of the $k$-tuple of primes falls within $\pgap(p)$.
Then in $\pgap(p)$ there is a constellation
$$G_1 = g_{20}\ldots g_{2n_2} g_{30}\ldots g_{3n_3} \ldots 
 g_{k0} \ldots g_{kn_k} $$
such that $\sum_{i=0}^{n_j} g_{ji} = b_j - b_{j-1}$.

As the recursion of \ref{RecursLemma} proceeds,
we identify a sequence of constellations $G_i$ which
ends with $G_N$ equal to the constellation $s_b$.
In the $\ord{i}$ iteration we pick an image of $G_i$
in which one addition $g_{jJ}+g_{j(J+1)}$ occurs.  
By Theorem~\ref{DelThm} this addition occurs, and
due to the size of $p$, this is the only addition that
occurs in this image of $G_i$; so the desired image
exists, and we take the resulting constellation as $G_{i+1}$.

For $N= n_2+n_3+\cdots+n_k$, $G_N$ is the constellation $s_b$.
\end{proof}

In \cite{Riesel}, Riesel uses Hardy and Littlewood's work \cite{HL}
to estimate the numbers of occurrences of the constellations
$424$ and $242$. Nicely \cite{NicelyQuads} has tested the estimate
for $242$ computationally.  The work in Section \ref{EstSection} 
supports the asymptotic order of the estimates.  In $\pgap(5)$ we
see the reason for Nicely's modularity condition $30n+11$.  By applying
the recursion a few times, we can sort copies of $242$ into branches
with more restrictive modularity conditions.  For any other constellation
$s$ we can observe analogous modularity conditions based on the location
of copies of $s$ in $\pgap(p)$.  Finally, $\pgap(5)$ supports the
conjecture \cite{Riesel} that $424$ should occur approximately twice
as often as $242$.

If our weaker conjecture \ref{surviveconj} is correct, we can provide solutions
to some problems posed by Erd\"os and Tur\'an \cite{ET}:
\begin{enumerate}
\item {\it Spikes.} What are $\limsup g_n/g_{n+1}$ and $\liminf g_n/g_{n+1}$?
\item {\it Oscillation.} Is there an $n_0$ such that for all $k \ge 1$,
 $g_{n_0+2k-1} < g_{n_0+2k}$ and $g_{n_0+2k} > g_{n_0+2k+1}$?
\item {\it Superlinearity.} Can $g_j < g_{j+1} < \ldots < g_{j+k}$ have
infinitely many solutions for every $k$?
\end{enumerate}

{\it Spikes.} We will construct a sequence $\{s_k\}$ of constellations, such
that $s_k$ occurs in $\pgap(p_k)$, $s_k$ is of the form $g_{k,n_k}2$,
and $g_{k,n_k}$ is strictly increasing in $k$.
Consider constellations in $\pgap(p_k)$ of the form
$g2$. For $s_k$ we take a pair with the maximal $g$.
In $\pgap(7)$, for example, this pair is $10,2$.  By Theorem~\ref{DelThm},
in $\pgap(p_{k+1})$ one of the images of this pair $g2$ has the gap $g$
added to the preceding gap.  This shows that the sequence $g_{k,n_k}$ is
strictly increasing in $k$.  Thus 
$$\limsup g_n/g_{n+1} \ge \limsup (g_{k,n_k}/2) = \infty.$$
By the symmetry of $\pgap(p)$,
$$\liminf g_n/g_{n+1} \le \liminf (2/g_{k,n_k}) = 0.$$

{\it Oscillation.} An affirmative answer to this question
asserts that eventually the gaps between primes oscillate in size.
We provide a counterexample.  In $\pgap(7)$, the constellation
$24682$ occurs.  Our conjecture implies that this constellation
will occur infinitely often as a constellation among prime
numbers.  This constellation is incompatible with the proposed
oscillation.

{\it Superlinearity.} This third problem inquires about
sequences of consecutive primes that exhibit superlinear growth.  
We provide
examples of the requested constellations, from which our
conjecture implies an affirmative answer to this problem.
From Remark~\ref{EasyRmk} we recall that in the middle of
$\pgap(p)$ there occurs a constellation of powers of $2$:
$2^k,\ldots,42424,\ldots,2^k$.  We can make this constellation
arbitrarily long by taking $p$ large enough.  The right half
of this constellation is an increasing sequence of $k$ gaps.
Assuming our conjecture holds, this constellation will occur
infinitely often as a constellation among primes.

\vskip .125in

Our work above also provides insight into the search for dense 
clusters of primes \cite{Riesel}.  One line of research looks for
constellations among the small primes to show up again later, among
larger primes.  Theorem~\ref{DelThm} and Lemma~\ref{SurviveLemma}
indicate limits on these searches.

\section{Conclusion}

We have studied directly the cycle of gaps
produced by each stage of Eratosthenes Sieve.
The work above can be divided into two parts: a deterministic
recursion for $\pgap(p)$, followed by statistical estimates
under an assumption of approximately uniform distributions.


Based on the recursion for the $p$-sieves, we have
conjectured that all constellations, which occur in $\pgap(p)$
for some prime $p$ and the sum of whose gaps is less than $2p$,
tend toward a uniform distribution in later stages of the sieve.
From this conjecture, we can make estimates of the number of
occurrences of a constellation between $p$ and $p^2$ for the
new prime $p$ at each stage of the sieve; all constellations
which occur before $p^2$ actually occur as constellations between primes.

We posed two conjectures, either of which imply that every sufficiently small 
constellation in $\pgap(p)$ occurs
infinitely often as a constellation among prime numbers.

Our stronger conjecture asserts that under the recursion of Eratosthenes sieve,
the images of a constellation are distributed approximately
uniformly in the fundamental cycle.
This conjecture allows us to estimate how many copies of
a constellation $s$ occur in the interval $[p,p^2]$, which estimates compare
favorably to the results from our initial computer searches.
For single gaps or for constellations consisting only of $2$'s and $4$'s, 
other estimates are available \cite{HL, Brent, Riesel,NicelyQuads}, and these
estimates agree at first order to ours.

This conjecture \ref{uniconj} on uniformity provides many open problems.
In our estimates, we use the interval $[1,\Pi_k]$ as the sample space
for the conjectured uniform distribution.
Simple observations about the structure of $\pgap(p)$ allow us to
adjust the sample space.
These refinements are on the order of $p$ and so improve our 
estimates only for very small prime numbers.
Are there any refinements that will improve our estimates in the large?

While the statistics in these estimates may be refined,
with the weaker Conjecture \ref{surviveconj} 
we have addressed three problems
posed by Erd\"os and Tur\'an \cite{ET}.


\bibliographystyle{amsplain}

\providecommand{\bysame}{\leavevmode\hbox to3em{\hrulefill}\thinspace}
\providecommand{\MR}{\relax\ifhmode\unskip\space\fi MR }
\providecommand{\MRhref}[2]{%
  \href{http://www.ams.org/mathscinet-getitem?mr=#1}{#2}
}
\providecommand{\href}[2]{#2}

\end{document}